\documentclass[11pt]{article}
\usepackage{amsthm}
\usepackage{amscd}
\usepackage{latexsym}
\usepackage{graphicx}
\usepackage{amsmath}
\usepackage{amsfonts}
\usepackage{float}
\usepackage{amssymb}
\usepackage{enumerate}
\usepackage{multirow}

\newtheorem{theorem}{Theorem}  [section]
\newtheorem{lemma}[theorem]{Lemma}

\newtheorem{definition}[theorem]{Definition}
\newtheorem{proposition}[theorem]{Proposition}

\theoremstyle{definition}

\theoremstyle{remark}

\begin{document}
\begin{center} \Huge \bf
The maximal angle between $3 \times 3$ copositive matrices
\end{center}
\bigskip
\begin{center}
\textsc{Daniel Gourion\footnote{Laboratoire de Mathématiques d'Avignon, Avignon Université, France
(e-mail address: daniel.gourion@univ-avignon.fr).} } 
 \end{center}
\medskip
\begin{quote} {\bf Abstract}.
In 2010, Hiriart-Urruty and Seeger posed the problem of finding the maximal possible angle $\theta_n$ between two copositive matrices of order $n$. They proved that $\theta_2=\frac{3}{4}\pi$. In this paper, we study the maximal angle between two copositive matrices of order 3. We show that $\theta_3=\frac{3}{4}\pi$ and give all possible pairs of matrices achieving this maximal angle. The proof is based on case analysis and uses optimization and basic linear algebra techniques.
\end{quote}
\begin{quote}{\bf Mathematics Subject Classification.}  15A48, 90C26. \end{quote}
\begin{quote}{\bf Key words.}  Convex cone, Copositive matrices,
Maximal Angle.\end{quote}

\section{Introduction}

We consider $\mathcal{S}_n$ the set of real symmetric matrices with the standard inner product
$$\langle A,B\rangle = Tr AB$$ and the associated norm $\|A\|= \sqrt{Tr A^2}$. A real symmetric matrix $A$ is copositive if $x^T A x\geq 0$ for every vector $x$ nonnegative entrywise. Let $\mathcal{P}_n$ be the set of positive semidefinite matrices, $\mathcal{N}_n$ be the set of entrywise nonnegative symmetric matrices and $\mathcal{C}_n$ the set of copositive matrices. These three sets are closed convex cones of $\mathcal{S}_n$. For any $n\geq 2$,  $\mathcal{P}_n+\mathcal{N}_n\subset \mathcal{C}_n$. Moreover it is known that for $n\leq 4$, $\mathcal{C}_n=\mathcal{P}_n+\mathcal{N}_n$ and that it is not true for $n\geq 5$ \cite{Di}. The cone of copositive matrices plays an important role in optimization (see for example \cite{Bo}, \cite{Du}) and its geometric structure has recently been studied in \cite {Af} and \cite{Dic}.

Let $\mathbb{S}_n$ denote the unit sphere of $\mathcal{S}_n$. The optimization problem
\begin{equation}\label{thetamax}
\theta_n:= \max_{A,B\in \mathcal{C}_n \cap\mathbb{S}_n}\arccos\langle A,B\rangle
\end{equation}
 defines the maximal angle of $\mathcal{C}_n$.  
 Interesting applications of the concept of maximal angle  of a closed convex cone arise in many areas of
mathematics, see \cite{ClLe,PeRe}. Ease computable formulas for maximal angle are given in Gourion and Seeger\,\cite[Proposition\,2]{GoSe}  for  Schur cones, in Seeger\,\cite[Example\,5.8]{Se}) for $p$\,-norm cones, in  Iusem and Seeger\,\cite[Theorem\,1]{IuSe} for ellipsoidal cones and in Gourion and Seeger \cite{GoSe2} for equilateral cones.

Computing the maximal angle of the set of copositive matrices is a challenging problem. Hiriart-Urruty and Seeger \cite{HiSe} showed that $\theta_2=3\pi/4$ and conjectured that $\theta_n= 3\pi/4$ for all $n\geq 2$. Goldberg and Shaked-Monderer \cite{GoSh} invalidated this conjecture, showing that $\lim_{n\rightarrow\infty}\theta_n=\pi$. Moreover, they computed the maximal angle $$\psi_n:=\max_{\substack{\scriptstyle A\in\mathcal{P}_n\cap\mathbb{S}_n   \\ B\in \mathcal{N}_n\cap\mathbb{S}_n}}\arccos\langle A,B\rangle$$ between a nonnegative matrix and a semidefinite positive matrix for $n=3,4$, proving that $\gamma_3=\gamma_4=3\pi/4$. Zhang \cite{Zh} showed that $\gamma_5=\arccos(-\frac{1+1/\sqrt{5}}{2})$. 
To our knowledge, the maximal angle $\theta_n$ between two copositive matrices for $n\geq 3$ is unknown, as well as the maximal angle $\gamma_n$ between $\mathcal{P}_n$ and $\mathcal{N}_n$ for $n\geq 6$.

In the sequel we will show that the maximal angle between two matrices in $\mathcal{C}_3$ is $3\pi /4$. The paper is structured as follows: in section 2, we recall some results concerning copositive matrices, in particular those of order 3 and we introduce some preliminary results about pairs $(A,B)$ of copositive matrices realizing a maximal angle. In section 3, we prove that $\theta_3=3\pi/4$ by case analysis. In section 4, we give our concluding remarks. 

\section{Preliminaries}
\subsection{Copositive matrices}

We first give some general properties of copositive matrices:

\begin{proposition}\label{neg} Let $A=(a_{ij})_{1\leq 1,j\leq n}$ be a copositive matrix of size $n$ then we have:
\begin{itemize}
\item[(i)] Every principal submatrix of A must also be copositive.
\item[(ii)] $a_{ii}\geq 0$ for all $1\leq i \leq n$.
\item[(iii)] $a_{ij}\geq -\sqrt{a_{ii}a_{jj}}$ for all $1\leq i < j \leq n$.
\item[(iv)] If $a_{ii}=0$ for some $1\leq i\leq n$, then $a_{ij}\geq 0$ for all $1\leq j \leq n$.
\item[(v)] If $a_{ij}<0$ for some $1\leq i<j\leq n$, then $a_{ii}> 0$ and $a_{jj}> 0$.
\item[(vi)] Suppose that the off-diagonal entries of $A$ are all nonpositive. Then $A\in\mathcal{P}_n$.
\item[(vii)] For any integers $i$ and $j$ such that $1\leq i \neq j \leq n$,
$-\frac{1}{2}\|A\|\leq a_{ij}\leq \frac{\sqrt{2}}{2}\|A\|.$
\end{itemize}
\end{proposition}

\begin{proof}
A proof for properties (i)-(v) can be found for example in \cite{Bu}, and (vi) comes from \cite{HiSe}. We prove (vii): from $a_{ij}^2+a_{ji}^2\leq \|A\|^2$ we deduce $a_{ij}\leq \sqrt{2}/2\times\|A\|$. On the other hand, any off-diagonal entry $a_{ij}$ of a copositive matrix $A$ is such that $a_{ij}+\sqrt{a_{ii}a_{jj}}\geq 0$, which, together with $2a_{ii}a_{jj}\leq a_{ii}^2+a_{jj}^2$, imply that any  entry is such that $a_{ij}\geq -1/2\times\|A\|$. 
\end{proof}

\noindent We now focus on $n=3$. Let $A=(a_{ij})_{1\leq 1,j\leq 3}$
 be a matrix in $\mathcal{S}_3$. Checking copositivity of such a matrix is easy to handle, using the following proposition, which can be found, for example in \cite{ChSe,Ha}.

\begin{proposition}
A symmetric matrix $A$ of order 3 is copositive if and only if the seven inequalities 
\begin{equation}\label{C1}
a_{11}\geq 0, a_{22}\geq 0, a_{33}\geq 0,
\end{equation}
\begin{equation}\label{C21}
\bar{a}_{12}:=a_{12}+\sqrt{a_{11}a_{22}}\geq 0,
\end{equation}
\begin{equation}\label{C22}
\bar{a}_{13}:=a_{13}+\sqrt{a_{11}a_{33}}\geq 0,
\end{equation}
\begin{equation}\label{C23}
\bar{a}_{23}:=a_{23}+\sqrt{a_{22}a_{33}}\geq 0,
\end{equation}
\begin{equation}\label{C3}
\sqrt{a_{11}a_{22}a_{33}}+a_{12}\sqrt{a_{33}}+a_{13}\sqrt{a_{22}}+a_{23}\sqrt{a_{11}}+\sqrt{2\bar{a}_{12}\bar{a}_{13}\bar{a}_{23}}\geq 0.
\end{equation}
are satisfied. 
\end{proposition}

\noindent The last inequality can be written in the disjunctive form 
$$
det A\geq 0 \quad \textrm{or} \quad \sqrt{a_{11}a_{22}a_{33}}+a_{12}\sqrt{a_{33}}+a_{13}\sqrt{a_{22}}+a_{23}\sqrt{a_{11}}\geq 0.$$

Now assume that all the diagonal elements of a symmetric matrix $A$ of order 3 are positive. Note that this is always true for matrices in $\mathcal{C}_3$ with at least 2 negative terms above the diagonal (see Prop. \ref{neg}). By diagonal scaling we introduce the auxiliary variables $\alpha,\beta$ and  $\gamma$ defined by $a_{12}=\alpha \sqrt{ a_{11}a_{22}}$, $a_{13}=\beta \sqrt{ a_{11}a_{33}}$ and $a_{23}=\gamma \sqrt{ a_{22}a_{33}}$. 
Inequalities  (\ref{C21}), (\ref{C22}), (\ref{C23}) and (\ref{C3}) are then respectively equivalent to 
\begin{equation}\label{C21'}
\alpha\geq -1,
\tag{3'}
\end{equation}
\begin{equation}\label{C22'}
\beta\geq -1,
\tag{4'}
\end{equation}
\begin{equation}\label{C23'}
\gamma\geq -1,
\tag{5'}
\end{equation}
\begin{equation}\label{C3'}
1+\alpha+\beta+\gamma+\sqrt{2(1+\alpha)(1+\beta)(1+\gamma)}\geq 0,
\tag{6'}
\end{equation}
with (6') which can be written in the disjunctive form 
$$
1+2\alpha\beta\gamma-\alpha^2-\beta^2-\gamma^2=(1-\alpha^2)(1-\beta^2) -(\gamma-\alpha\beta)^2 \geq 0   \quad \textrm{or} \quad 1+\alpha+\beta+\gamma\geq 0,$$ the first inequality being equivalent to $
\det A\geq 0$.

 Let $A\in\mathcal{S}_3$ and $a\in\{0,1,2,3\}$ be the number of negative elements of $A$ above the diagonal. In the following, without loss of generality, we will assume that when $a=1$ then $a_{12}<0$ (that is $\alpha<0$) and when $a=2$, $a_{12}<0$ and  $a_{13}<0$ (that is $\alpha<0$ and $\beta<0$). We now give necessary and sufficient conditions for a matrix $A$ to be in $\mathcal{C}_3$ according to its number of negative elements.

\begin{proposition}\label{cascopsign}
Let $A\in\mathcal{S}_3$ with non-negative elements on its diagonal and $a\in\{0,1,2,3\}$ be the number of negative elements of $A$ above its diagonal. Then
\begin{itemize}
\item If $a=0$ then $A\in\mathcal{C}_3$.
\item If $a=1$, then $A\in\mathcal{C}_3$ if and only if $a_{12}\geq -\sqrt{a_{11}a_{22}}$.
\item If $a=2$, then $A\in\mathcal{C}_3$ if and only if $\alpha\geq -1$, $\beta\geq -1$ and $\gamma\geq \alpha\beta-\sqrt{(1-\alpha^2)(1-\beta^2)}.$ Moreover, in this case, $A\in\mathcal{P}_3$ if and only if $\det A\geq 0$, $\alpha\geq -1$ and $\beta\geq -1$.
\item If $a=3$,  then $A\in\mathcal{C}_3$ if and only if $A\in\mathcal{P}_3$.
\end{itemize}
\end{proposition}

\begin{proof}
Let $A\in\mathcal{S}_3$ with non-negative elements on its diagonal. The only nontrivial case is when $a=2$. In this case, first suppose that $\alpha\geq -1$, $\beta\geq -1$ and $\gamma\geq \alpha\beta-\sqrt{(1-\alpha^2)(1-\beta^2)}$. Then inequalities (\ref{C1}), (\ref{C21'}), (\ref{C22'}) and (\ref{C23'}) are verified. Moreover, we either have $|\gamma-\alpha\beta| \leq \sqrt{(1-\alpha^2)(1-\beta^2)}$ or $\gamma> \alpha\beta +\sqrt{(1-\alpha^2)(1-\beta^2)}$. In the first case we immediately deduce that $\det A \geq 0$ and in the second case, we deduce that $$1+\alpha+\beta+\gamma > 1+\alpha+\beta +\alpha\beta =(1+\alpha)(1+\beta) \geq 0 .$$ Hence, in both cases, inequality (\ref{C3'}) is verified and $A$ is copositive. We now prove the converse statement. Suppose that $A$ is copositive with $\alpha<0$, $\beta<0$, $\gamma\geq 0$. Inequalities (\ref{C21'}), (\ref{C22'}) imply  $\alpha\geq -1$, $\beta\geq -1$. Now suppose that $\gamma<\alpha\beta-\sqrt{(1-\alpha^2)(1-\beta^2)}$. Then $\det A<0$ and
$$1+\alpha+\beta+\gamma<1+\alpha+\beta+\alpha\beta-\sqrt{(1-\alpha^2)(1-\beta^2)}.$$ We now observe that for any $-1\leq\alpha,\beta<0$, we have $$(1+\alpha)(1+\beta)\leq\sqrt{(1-\alpha^2)(1-\beta^2)}.$$ Thus $1+\alpha+\beta+\gamma <0$ and $A$ is not copositive. We now prove the second assertion. Suppose $\det A \geq 0$. We are going to prove that all principal minors of $A$ are nonnegative. Because $\alpha \geq -1$ and $\beta \geq -1$, we only have to prove that $\gamma\leq 1$. On one hand, $\det A\geq 0$ implies that $\gamma\leq \alpha\beta+\sqrt{(1-\alpha^2)(1-\beta^2)}$. On the other hand, observe that $(1-\alpha^2)(1-\beta^2)\leq(1-\alpha\beta)^2$. We immediately deduce $\gamma\leq\alpha\beta+\sqrt{(1-\alpha^2)(1-\beta^2)}\leq 1$. Hence $A$ is positive semidefinite. The converse assertion is trivial. This completes the proof for $a=2$.

\noindent When $a=0$, then $A\in\mathcal{N}_3\subset\mathcal{C}_3$. When $a=1$, $a_{12}\geq -\sqrt{a_{11}a_{22}}$, then conditions (2) to (6) are immediately checked. When $a=3$ and $A\in\mathcal{C}_3$, then Prop. (2.1, vi) implies that $A\in\mathcal{P}_3$.
\end{proof}

\subsection{Maximal angle}

We now introduce the notion of antipodal pair, as in \cite{IuSe2}. 

\begin{definition}
 $A^\star$ and $B^\star$ are said to be antipodal in $\mathcal{C}_n$ if $A^\star,B^\star\in\mathcal{C}_n\cap\mathbb{S}_n$ and $$\arccos\langle A^\star,B^\star\rangle
=\max_{A,B\in \mathcal{C}_n \cap\mathbb{S}_n}\arccos\langle A,B\rangle.$$

\end{definition}

\noindent Please note that because $\arccos$ is a decreasing function on $[-1,1]$, maximizing the angle $\arccos \langle A,B\rangle$ is the same as minimizing the inner product $\langle A,B\rangle$. The following proposition gives a condition on the signs of the entries of an antipodal pair of $\mathcal{C}_n$. This condition will help us to build a proof by case analysis to determine the maximal angle of $\mathcal{C}_3$. 
\begin{proposition}\label{sign}
Let $(A,B)$ be an antipodal pair in $\mathcal{C}_n$. Then the following results hold for any integers $i$ and $j$ such that $1\leq i \neq j \leq n$:
\begin{itemize}
\item[(i)]If $a_{ij}<0$ then $b_{ij}>0$.
\item[(ii)]If $a_{ij}=0$ then $b_{ij}\geq 0$.
\end{itemize}
\end{proposition}

\begin{proof}
Let  $A,B \in \mathcal{C}_n\cap\mathbb{S}_n$ and $i,j$ be two integers  such that $1\leq i \neq j \leq n$ and  $a_{ij}<0$, $b_{ij}\leq 0$. Take $A'$ with the same entries as $A$ except $a'_{ij}=a'_{ji}=0$. Then $(1/\|A'\|)A'\in\mathcal{C}_n+\mathcal{N}_n = \mathcal{C}_n$ and the angle between $(1/\|A'\|)A'$ and $B$ is greater than the angle between $A$ and $B$. That proves the first assertion. The contraposition of this assertion immediately induces the second assertion.
\end{proof}

\noindent The following results are well known and can be found in \cite{ChSe,Ha}:

\begin{proposition}
$\theta_2=\psi_2=3\pi/4$. Furthermore, the pair 
\begin{eqnarray*}
 A =\frac{1}{2}\left[\begin{array}{cc}
1 & -1 \\
-1& 1\\
\end{array}\right]
, \quad B = \frac{1}{2} \left[\begin{array}{cc}
0 &\sqrt{2}\\
\sqrt{2} & 0\\
\end{array}\right]\,
\end{eqnarray*} is the only antipodal pair in $\mathcal{C}_2$.
\end{proposition}

\begin{lemma}
$\theta_3\geq 3\pi/4$.
\end{lemma}

\noindent This lemma is obvious because $(\theta_n)_n$ is an increasing sequence. For example, the angle between the following pair $(A,B)$ is equal to $3\pi/4$.
\begin{eqnarray*}
 A =\frac{1}{2}\left[\begin{array}{ccc}
1 & -1 & 0 \\
-1& 1 & 0 \\
0& 0 & 0 \\
\end{array}\right]
, \quad B = \frac{1}{2}\left[\begin{array}{ccc}
0 &\sqrt{2} & 0\\
\sqrt{2} & 0 & 0\\
0& 0 & 0 \\
\end{array}\right]\,
\end{eqnarray*}

\section{The maximal angle between copositive matrices of size 3} 

In this section we will investigate necessary conditions verified by an antipodal pair $(A,B)$  with respectively $a$ and $b$ negative entries above the main diagonal. This will lead us to the following conclusion:

\begin{theorem}
The maximal angle between two copositive matrices of size 3 is equal to $3\pi/4$. This value is attained only (up to row and column permutation and multiplication by a positive scalar) by pairs defined by: 

\begin{eqnarray}
 A =\left[\begin{array}{ccc}
1/2 & -\sqrt{a_{22}/2} & -\sqrt{a_{33}/2} \\
-\sqrt{a_{22}/2} & a_{22} & \sqrt{a_{22}a_{33}}\\
-\sqrt{a_{33}/2} & \sqrt{a_{22}a_{33}} & a_{33}\\
\end{array}\right]
, \quad B = \left[\begin{array}{ccc}
0 & \sqrt{a_{22}} & \sqrt{a_{33}} \\
\sqrt{a_{22}}& 0 & 0\\
\sqrt{a_{33}}  & 0 & 0\\
\end{array}\right],\,
\end{eqnarray}
with $a_{22}\geq 0$, $a_{33}\geq 0$ and $a_{22}+a_{33}=1/2$.
\end{theorem} 

The proof is just a matter of gathering the results of the following subsections. Without loss of generality, we suppose that $a\geq b$. Prop. \ref{sign} induces that $a+b\leq 3$. Hence $(a,b)$ can take only 6 values, that are $(0,0)$, $(1,0)$, $(1,1)$, $(2,0)$, $(2,1)$ and $(3,0)$. Two of these six cases, namely $(0,0)$, $(1,0)$  are easy to handle. Cases $(2,0)$ and $(3,0)$ are solved by linear algebra techniques. Finally, the last two cases $(1,1)$ and $(2,1)$ are handled by optimization techniques. As the calculations in these last two subsections are fairly complex, they have been verified using both the computer algebra system Xcas and the numerical software Matlab.

\subsection{Case $(0,0)$}

Here both matrices being entrywise non-negative, their inner product is non-negative and we immediately deduce the following proposition:
\begin{proposition} For any couple $(A,B)$ of copositive matrices entrywise non-negative, the angle between $A$ and $B$ is less than or equal to $\pi/2$.
\end{proposition}

\subsection{Case $(1,0)$}

\begin{proposition}\label{case10}
For any couple $(A,B)$ of copositive matrices such that $A$ has exactly one negative element above its diagonal and $B$ is non-negative entrywise, the angle between $A$ and $B$ is not greater than $3\pi/4$ and this upper bound is attained only by this couple (up to row and column permutation and multiplication by a positive scalar): 

\begin{eqnarray}\label{solsimple}
 A =\frac{1}{2}\left[\begin{array}{ccc}
1 & -1 & 0 \\
-1& 1 & 0 \\
0& 0 & 0 \\
\end{array}\right]
, \quad B =\frac{1}{2} \left[\begin{array}{ccc}
0 &\sqrt{2} & 0\\
\sqrt{2} & 0 & 0\\
0& 0 & 0 \\
\end{array}\right]\,
\end{eqnarray}
\end{proposition}

\begin{proof}
It is directly derived from Prop. \ref{neg} (vii).
\end{proof}
\subsection{Case $(3,0)$}

\noindent Prop. \ref{neg} implies that if  $(a,b)=(3,0)$, then $A\in\mathcal{P}_3$ and $B\in\mathcal{N}_3$. We deduce immediately that in this case the angle between $A$ and $B$ is not greater than $\psi_3=3\pi/4$. In fact, we can even prove that this upper bound is not attained:

\begin{proposition} For any couple $(A,B)$ of copositive matrices with all off-diagonal entries of $A$ negative and $B$ entrywise non-negative, the angle between $A$ and $B$ is less than $3\pi/4$. Moreover, the supremum of this set of angles is equal to $3\pi/4$.
\end{proposition}

\begin{proof}
Let $A$ and $B$ be two copositive matrices such that $a=3$, $b=0$ and $\|A\|=\|B\|=1$.  We have $B\in\mathcal{N}_3$ and from Prop 2.3 we know that $A\in\mathcal{P}_3$.
Defining 
\begin{eqnarray*}
 B_A =\frac{1}{\sqrt{2(a_{12}^2+a_{13}^2+a_{23}^2)}}\left[\begin{array}{ccc}
0 & -a_{12} & -a_{13} \\
-a_{12} & 0 & -a_{23}\\
-a_{13} & -a_{23} & 0\\
\end{array}\right]\,
\end{eqnarray*}
\\
\noindent and using  Cauchy–Schwarz inequality, we obtain $\langle A,B\rangle\geq \langle A,B_A\rangle=-\sqrt{2(a_{12}^2+a_{13}^2+a_{23}^2)}$. Note that $B_A\in\mathcal{N}_3\setminus\mathcal{P}_3$ and $\|B_A\|=1$. Denote by $\Lambda$ and $\Lambda_-$ respectively the sets of eigenvalues (resp. negative eigenvalues) of $B_A$. Following Prop 2.1 in \cite{GoSh}, we deduce that $$\langle A,B_A\rangle \geq -\frac{\sqrt{\sum_{\lambda\in\Lambda_-} \lambda^2}}{\sqrt{{\sum_{\lambda\in\Lambda} \lambda^2}}}.$$ Observe that $\det B_A=-2a_{12}a_{13}a_{23}>0$ and tr$(B_A)=0$ and we infer that $B_A$ has one positive eigenvalue $\lambda_1$ and two negative eigenvalues $\lambda_2, \lambda_3$ such that $\lambda_1=-(\lambda_2+\lambda_3)$. This implies that $$\langle A,B_A\rangle \geq -\frac{\sqrt{\lambda_2^2+\lambda_3^2}}{\sqrt{2\lambda_2^2+2\lambda_3^2+2\lambda_2\lambda_3}}>-1/\sqrt{2}.$$ Hence the angle between $A$ and $B$ is less than $3\pi/4$. We prove that this upper bound is a supremum with a continuity argument between the case (3,0) and the case (1,0):  the angle between copositive matrices $A_\epsilon$ and $B$ defined by 
\begin{eqnarray*}
 A_\epsilon =\left[\begin{array}{ccc}
1 & -1+3\epsilon & -\epsilon \\
-1+3\epsilon & 1 & -\epsilon\\
-\epsilon & -\epsilon & \epsilon\\
\end{array}\right]
, \quad B =\frac{1}{2} \left[\begin{array}{ccc}
0 &\sqrt{2} & 0\\
\sqrt{2} & 0 & 0\\
0& 0 & 0 \\
\end{array}\right],\quad  0\leq\epsilon < 1/3
\end{eqnarray*}
tends to $3\pi/4$ when $\epsilon$ tends to 0 from above. The condition $0\leq\epsilon < 1/3$ is given so that $A_\epsilon$ is copositive with negative terms above its diagonal.
\end{proof}

\subsection{Case $(2,0)$}

This case is a little bit more tedious: the study of the maximal angle differs according to whether $A$ is positive semidefinite or not. Without loss of generality, suppose $a_{12}< 0$, $a_{13}< 0$, $a_{23}\geq 0$. We can state the following proposition: 

\begin{proposition}\label{case20}
For any couple $(A,B)$ of copositive matrices with respectively two and no negative element above the diagonal, the angle between $A$ and $B$ is not greater than $3\pi/4$ and this upper bound is attained only (up to row and column permutation and multiplication by a positive scalar) by pairs defined by: 

\begin{eqnarray}
 A =\left[\begin{array}{ccc}
1/2 & -\sqrt{a_{22}/2} & -\sqrt{a_{33}/2} \\
-\sqrt{a_{22}/2} & a_{22} & \sqrt{a_{22}a_{33}}\\
-\sqrt{a_{33}/2} & \sqrt{a_{22}a_{33}} & a_{33}\\
\end{array}\right]
, \quad B = \left[\begin{array}{ccc}
0 & \sqrt{a_{22}} & \sqrt{a_{33}} \\
\sqrt{a_{22}}& 0 & 0\\
\sqrt{a_{33}}  & 0 & 0\\
\end{array}\right],\,
\end{eqnarray}
with $a_{22}>0$, $a_{33}>0$ and $a_{22}+a_{33}=1/2$.
\end{proposition}

\begin{proof}
Let $A$ and $B$ be two copositive matrices such that $a=2$, $b=0$ and $\|A\|=\|B\|=1$. Without loss of generality, take $a_{12}< 0$, $a_{13}< 0$, $a_{23}\geq 0$.  Using  Cauchy–Schwarz inequality and defining 
\begin{eqnarray*}
 B_A =\frac{1}{\sqrt{2(a_{12}^2+a_{13}^2)}}\left[\begin{array}{ccc}
0 & -a_{12} & -a_{13} \\
-a_{12} & 0 & 0\\
-a_{13} & 0 & 0\\
\end{array}\right]\,
\end{eqnarray*}
\\
\noindent we get $\langle A,B\rangle\geq \langle A,B_A\rangle=-\sqrt{2(a_{12}^2+a_{13}^2)}$ for all $B$ copositive such that $b=0$ and $\|B\|=1$. Note that $B_A\in\mathcal{N}_3\setminus\mathcal{P}_3$ and $\|B_A\|=1$.

\noindent We first suppose that $A\in\mathcal{P}_3$. Using the same logic as in the previous subsection, we obtain that $B_A$ has three eigenvalues, that are $-\lambda_1, 0, \lambda_1$ with $\lambda_1>0$. Hence $$\langle A,B_A\rangle \geq -\frac{\sqrt{\lambda_1^2}}{\sqrt{2\lambda_1^2}}=\frac{-1}{\sqrt{2}},$$ with equality only when $A$ is the negative definite part of $B_A$, up to multiplication by a positive scalar (see Prop 2.1 in \cite{GoSh}). Computing this negative definite part and normalizing it, we obtain that $\langle A,B_A\rangle =\frac{-1}{\sqrt{2}}$ if and only if $A$ verifies the following equation:
\begin{eqnarray*}
 A =\frac{1}{2}\left[\begin{array}{ccc}
1 & a_{12}/\sqrt{(a_{12}^2+a_{13}^2)} & a_{13}/\sqrt{(a_{12}^2+a_{13}^2)} \\
a_{12}/\sqrt{(a_{12}^2+a_{13}^2)} & a_{12}^2/(a_{12}^2+a_{13}^2) & a_{12}a_{13}/(a_{12}^2+a_{13}^2)\\
a_{13}/\sqrt{(a_{12}^2+a_{13}^2)} & a_{12}a_{13}/(a_{12}^2+a_{13}^2) & a_{13}^2/(a_{12}^2+a_{13}^2)\\
\end{array}\right]\,.
\end{eqnarray*}
Solving this equation gives $2\sqrt{a_{12}^2+a_{13}^2}=1$, $a_{11}=1/2$, $a_{22}=2a_{12}^2$, $a_{33}=2a_{13}^2$ and $a_{23}=2a_{12}a_{13}$, which is equivalent to the expression of $A$ given in (8).

\noindent We still have to deal with the case where $A\notin\mathcal{P}_3$. $A$ being copositive, we infer from Prop. 2.3 that $\det A<0$. Using the auxiliary variables $\alpha,\beta,\gamma$ and Prop. 2.3 we know that it is the case if and only if
$$
(\gamma-\alpha\beta)^2>(1-\alpha^2)(1-\beta^2).
$$
From Prop. \ref{cascopsign} we deduce that $\gamma>\alpha\beta+\sqrt{(1-\alpha^2)(1-\beta^2)}$. Now take $\gamma'$ such that  $\gamma'=\alpha\beta+\sqrt{(1-\alpha^2)(1-\beta^2)}$ and define a matrix $A'$ with parameters $a_{11},a_{22},a_{33},\alpha, \beta, \gamma'$ and take $\tilde{A}=\frac{1}{\|A'\|}A'$. These matrices are copositive and such that $\|A'\|<1$,  $\tilde{a}_{12}<a'_{12}=a_{12}<0$, $\tilde{a}_{13}<a'_{13}=a_{13}<0$, $a_{23}>\tilde{a}_{23}>a'_{23}> 0$. Hence, $\tilde{A}$ has exactly two negative elements above the diagonal and $\det \tilde{A}=0$ which means that $\tilde{A}\in\mathcal{P}_3$.
Moreover we observe that $\langle A,B_A\rangle > \langle \tilde{A},B_{\tilde{A}}\rangle\geq \frac{-1}{\sqrt{2}}$. Thus, no antipodal pair can be found when $A\notin\mathcal{P}_3$.

\end{proof}

\subsection{Case $(1,1)$}
Here we are looking for an antipodal pair $(A,B)$ such that both $A$ and $B$ have exactly one negative entry above the diagonal. Without loss of generality, we suppose that $a_{12}< 0$, $a_{13}\geq 0$, $a_{23}\geq 0$, $b_{12}\geq 0$, $b_{13} < 0$, $b_{23}\geq 0$. Following proposition \ref{sign}, this directly implies, that if $(A,B)$ is an antipodal pair, then $b_{12}>0$ and $a_{13}>0$. Moreover, conditions (\ref{C21}) and (\ref{C22}) imply that $a_{11}$, $a_{22}$, $b_{11}$ and $b_{33}$ are positive. The following lemma gives some necessary conditions for such a pair $(A,B)$ to be an antipodal pair:

\begin{lemma}
Let $(A,B)$ be an antipodal pair as described above. Then $a_{33}=a_{23}=b_{22}=b_{23}=0$, $a_{12}=-\sqrt{a_{11}a_{22}}$ and $b_{13}=-\sqrt{b_{11}b_{33}}$.
\end{lemma}
\begin{proof}
If $a_{33}\neq 0$ then consider a matrix $A'$ equal to $A$ except $a'_{33}=0$. Then $A'$ is copositive (Prop. 2.3) and $\|A'\|<1$. Using the fact that the inner product is negative for an antipodal pair, we deduce that the angle between $A'$ and $B$ is  greater than the angle between $A$ and $B$. We prove similarly that $b_{22}=a_{23}=b_{23}=0$. 
Now we prove that $a_{12}=-\sqrt{a_{11}a_{22}}$. Suppose that $-\sqrt{a_{11}a_{22}}<a_{12}<0$. Then, for a sufficiently small $\epsilon>0$, define $\bar{A}$ with the same entries as $A$ except $\bar{a}_{11}=\sqrt{a_{11}^2-\epsilon}$
 and $\bar{a}_{12}^2=-\sqrt{a_{12}^2+\epsilon/2}$. Then 
the angle between $\bar{A}$ and $B$ is strictly greater that the angle between $A$ and $B$. We prove in a similar way that $b_{13}=-\sqrt{b_{11}b_{33}}$.
\end{proof}

\begin{proposition} Let $(A,B)$ be a couple of copositive matrices with each one containing exactly one negative  element above the diagonal. Then $(A,B)$ is not an antipodal pair for problem (1). Moreover, it is not a local maximum for problem (1).  Finally, the supremum of this set of angles is equal to $3\pi/4$.
\end{proposition}

\begin{proof}
\noindent Suppose that $(A,B)$ is an antipodal pair with $A$ and $B$ containing exactly one negative  element above the diagonal. We use  the previous lemma and the norm condition $\|A\|=\|B\|=1$ to obtain $(a_{11}+a_{22})^2+2a_{13}^2=1$ and $(b_{11}+b_{33})^2+2b_{12}^2=1$. $(A,B)$ being an antipodal pair, then it is an optimal solution of the following optimization problem:

\begin{equation}\label{opt11}
   \min_{\substack{\scriptstyle a_{11}>0, a_{22}>0, b_{11}>0, b_{33}>0 \\ a_{11}+a_{22}<1,b_{11}+b_{33}<1 }}\quad a_{11}b_{11}-\sqrt{2a_{11}a_{22}}\sqrt{1-(b_{11}+b_{33})^2}-\sqrt{2b_{11}b_{33}}\sqrt{1-(a_{11}+a_{22})^2},
\end{equation}
where the function to be minimized is the inner product between two matrices as described in the previous lemma.  We now look if there is a minimizer of the objective function in the admissible set of (\ref{opt11}).

\noindent For that purpose, we solve the first order necessary optimality conditions and we  obtain
$$
\frac{b_{11}}{b_{33}}=\frac{a_{11}}{a_{22}}=\frac{\sqrt{5}-1}{2},
$$
which leads to a unique solution
$$
a_{11}=b_{11}=\sqrt{\frac{7\sqrt{5}-15}{10}},\quad a_{22}=b_{33}=\sqrt{\frac{3\sqrt{5}-5}{10}} ,\quad a_{12}=b_{13}=-\sqrt{\frac{5-2\sqrt{5}}{5}},\quad a_{13}=b_{12}=\sqrt{\frac{\sqrt{5}-1}{2\sqrt{5}}}
$$
corresponding to a inner product $$\langle A,B\rangle =\frac{1-\sqrt{5}}{2}\approx -0.6180$$ and an angle between $A$ and $B$ equal to $\arccos\langle A,B\rangle \approx 0.7121\pi$. We conclude that there is no antipodal pair for problem (1) in case (1,1). Actually, we have checked that this pair is a critical pair but not a local minimum: the Hessian matrix of the objective function at this point admits a negative eigenvalue. 

\noindent It is easy to see that for this problem, the value of the objective function on the boundary of the admissible set relatively to $\mathcal{C}_3\cap\mathbb{S}_3$ is greater than or equal to $-\sqrt{2}/2$. Take for example $a_{11}=a_{22}=1/2,b_{11}=b_{33}=0$, which is the solution given in (\ref{solsimple}), for reaching this value. This proves the result about the supremum.

\end{proof}

\subsection{Case $(2,1)$}

Without loss of generality we will choose $a_{12}<0$, $a_{13}<0$ and $b_{23}<0$. Suppose that $(A,B)$ is antipodal. 
In this case, proposition \ref{sign} together with conditions \ref{C21}, \ref{C22} and \ref{C23} immediately imply that it verifies $a_{11}>0, a_{22}>0, a_{33}>0, a_{23}>0, b_{22}>0, b_{33}>0, b_{12}>0$ and $b_{13}>0$. More details on a possible antipodal pair can be found in the following lemma:

\begin{lemma}
Let $(A,B)$ be an antipodal pair as described above. Then $b_{11}=0$, $b_{23}=-\sqrt{b_{22}b_{33}}$ and $a_{23}>\sqrt{a_{22}a_{33}}$.
\end{lemma}

\begin{proof}
The two equalities can be proved using the same reasoning as in the previous section.  We now show that $a_{23}>\sqrt{a_{22}a_{33}}$. Suppose to the contrary that  $a_{23}\leq \sqrt{a_{22}a_{33}}$. We take $B'$ with the same entries as $B$ except that $b_{22}=b_{33}=b_{23}=0$. Then $B'\in\mathcal{C}_3$, $\|B'\|<1$ and 

\begin{align*}
   0  > \langle A,B \rangle & = a_{22}b_{22}+a_{33}b_{33}+2a_{12}b_{12}+2a_{13}b_{13}-2a_{23}\sqrt{b_{22}b_{33}}  \\
   & \geq 2a_{12}b_{12}+2a_{13}b_{13}+(\sqrt{a_{22}b_{22}}-\sqrt{a_{33}b_{33}})^2 \\
    & \geq 2a_{12}b_{12}+2a_{13}b_{13}=\langle A,B' \rangle \\
    & > \frac{1}{\|B'\|}\langle A,B' \rangle.
\end{align*}
\noindent Hence $(A,B)$ is not an antipodal pair and we have proved by contradiction that $\gamma>1$. 
\end{proof}

\begin{proposition} For any couple $(A,B)$ of copositive matrices with respectively two and one  negative  element above the diagonal, the angle between $A$ and $B$ is less than $3\pi/4$. Moreover, the supremum of this set of angles is equal to $3\pi/4$.
\end{proposition}

\begin{proof}

\noindent If $(A,B)$ is an antipodal pair with $a=2$ and $b=1$, it checks all the above conditions, and $B$ is an optimal solution  of the following problem:

\begin{equation}
   \min_{\substack{\scriptstyle B\in \mathbb{S}_3 \\ b_{22}>0,b_{33}>0,b_{12}>0,b_{13}>0}}\quad \langle A,B \rangle=a_{22}b_{22}+a_{33}b_{33}+2a_{12}b_{12}+2a_{13}b_{13}-2a_{23}\sqrt{b_{22}b_{33}}.
\end{equation}

\noindent Computing the first order optimality conditions, we obtain that $B$ has to verify:

\[
\left \{
\begin{array}{c @{=} c}
    a_{22}b_{13}\sqrt{b_{22}}-a_{13}(b_{22}+b_{33})\sqrt{b_{22}}-a_{23}b_{13}\sqrt{b_{33}} & 0 \\
    a_{33}b_{13}\sqrt{b_{33}}-a_{13}(b_{22}+b_{33})\sqrt{b_{33}}-a_{23}b_{13}\sqrt{b_{22}} & 0 \\
    a_{12}b_{13}-a_{13}b_{12} & 0 
\end{array}
\right.
\]

\noindent We solve these equations by introducing  $$x=\sqrt{\frac{b_{33}}{b_{22}}}.$$ We obtain a quadratic equation in $x$, the unique relevant solution of which is $$x=\frac{a_{22}-a_{33}+\sqrt{(a_{33}-a_{22})^2+4 a_{23}^2}}{2a_{23}}.$$ We deduce that there is a unique critical point which is as follows:
$$
   b_{12}  =  -\frac{a_{12}}{\sqrt{(a_{22}-a_{23}x)^2+2a_{12}^2+2a_{13}^2}},\quad    
  b_{13}  =  \frac{b_{12}a_{13}}{a_{12}} ,\quad b_{22}=\frac{b_{12}(a_{22}-a_{23}x)}{a_{12}(1+x^2)},\quad b_{33}=b_{22}x^2,
$$
and we will denote it by $B_A$. We have shown that if $(A,B)$ is antipodal then $B=B_A$. Then, we define $f$ by $f(A)=\langle A,B_A \rangle$. If $(A,B_A)$ is antipodal, then $A$ minimizes $f$ over the set $$E:=\{A\in\mathcal{C}_3\cap\mathbb{S}_3:a_{ii}>0\forall i\in\{1,2,3\},-\sqrt{a_{11}a_{22}}\leq a_{12}<0,-\sqrt{a_{11}a_{33}}\leq a_{13}<0, a_{23}>\sqrt{a_{22}a_{33}}\}.$$ We
 now reformulate $f$ and, after somewhat cumbersome computations, we obtain:
\begin{align*}
f(A) & =-\frac{1}{\sqrt{2}}\sqrt{a_{22}^2+a_{33}^2+4a_{12}^2+4a_{13}^2+2a_{23}^2-(a_{22}+a_{33})\sqrt{(a_{33}-a_{22})^2+4a_{23}^2}}\\ & :=  -\frac{1}{\sqrt{2}}\sqrt{h(A)}.
\end{align*}

\noindent The fact that $\|A\|=1$ implies that minimizing $f$ is equivalent to minimizing $g$ defined by
$$
g(A):=2a_{11}^2+a_{22}^2+a_{33}^2+2a_{23}^2+(a_{22}+a_{33})\sqrt{(a_{33}-a_{22})^2+4a_{23}^2}$$ 
We deduce that if this problem has an optimal value, then $a_{12}=-\sqrt{a_{11}a_{22}}$ and $a_{13}=-\sqrt{a_{11}a_{33}}$. In fact suppose that $a_{12}>-\sqrt{a_{11}a_{22}}$. Then we can build $A'$ by slightly decreasing $a_{11}$ and possibly $|a_{13}|$ and increasing $|a_{12}|$ in such a way as to keep the norm of $\|A'\|$ equal to 1 and to have $g(A')<g(A)$.

\noindent Thus, let $A$ be a matrix in $E$ verifying $a_{12}=-\sqrt{a_{11}a_{22}}$ and $a_{13}=-\sqrt{a_{11}a_{33}}$. Now, 
 observe that
$g(A)+h(A)=2$ and $$g(A)-h(A)=2a_{11}^2-4a_{11}(a_{22}+a_{33})+2(a_{22}+a_{33})\sqrt{(a_{33}-a_{22})^2+4a_{23}^2}.$$

The discriminant of the quadratic polynomial $g(A)-h(A)$ in the variable $a_{11}$ is negative. Hence $g(A)-h(A)>0$ for any $A$, resulting in $g(A)>1$, $h(A)<1$ and finally $f(A)>-1/\sqrt{2}$. We deduce that $\arccos\langle A,B\rangle< 3\pi/4$, for any antipodal pair $(A,B)$ such that $a=2$ and $b=1$, which is a contradiction. Hence, there is no antipodal pair in the case (2,1). The result about the supremum is left  to the reader, as it follows the same lines as in case (3,0).

%

For the sake of completeness, we conclude this case study by mentioning that the set of pairs $(A,B)$ verifying  the first order optimality conditions for minimizing $g$ can be computed here as in case (1,1), leading to quite complicated expressions.
\end{proof}

\section{Concluding remark}

We have shown by case analysis that the maximal angle between copositive matrices of order three is equal to $3\pi/4$ and found all the pairs $(A,B)$ of copositive matrices such that $\langle A,B\rangle =\cos 3\pi/4 = -\sqrt{2}/2$. Unfortunately, the nature of the proof makes it impossible to generalize to any value of $n$. Even the case $n=4$ could prove very tedious, but maybe possible to handle: a similar approach by case analysis would lead to no less than 16 cases to analyze, with not much of them being easy.

\vspace{0.5cm}

\noindent \textbf{Conflict of interest statement.}
The author Daniel Gourion certifies that he has no affiliations with or involvement in any organization or entity with any financial interest (such as honoraria; educational grants; participation in speakers’ bureaus;
membership, employment, consultancies, stock ownership, or other equity interest; and expert testimony or patent-licensing
arrangements), or non-financial interest (such as personal or professional relationships, affiliations, knowledge or beliefs) in
the subject matter or materials discussed in this manuscript.

\vspace{0.5cm}

\noindent\textbf{Data availability statement.}
The author confirms that the data supporting the findings of this study are available within the article.

\end{document}